\documentclass[12pt,thmsa]{article}%
\usepackage{amsmath}
\usepackage{amssymb}
\usepackage{sw20bams}%
\usepackage{amsfonts}%
\usepackage{graphicx}
\providecommand{\U}[1]{\protect\rule{.1in}{.1in}}
\begin{document}

\author{Steven Finch and Pascal Sebah}
\title{Squares and Cubes Modulo $n$}
\date{March 25, 2016}
\maketitle

\begin{abstract}
We study the asymptotics of the average number of squares (or quadratic
residues)\ in $\mathbb{Z}_{n}$ and $\mathbb{Z}_{n}^{*}$. Similar analyses are
performed for cubes, square roots of $0$ and $1$, and cube roots of $0$ and
$1$.

\end{abstract}

\footnotetext{Copyright \copyright \ 2006 by Steven R. Finch. All rights
reserved.}Let $\mathbb{Z}_{n}$ denote the ring of integers modulo $n$, and let
$\mathbb{Z}_{n}^{*}$ denote the group (under multiplication) of integers
relatively prime to $n$. The number of elements in $\mathbb{Z}_{n}^{*}$ is
$\varphi(n)$, where $\varphi$ is Euler's totient function. What is the average
number of elements in $\mathbb{Z}_{n}^{*}$, given an arbitrary $n$? One way to
answer this question is to apply the Selberg-Delange method \cite{Slbrg,
Dlnge, Tnbm} to the Dirichlet series
\begin{align*}%
{\displaystyle\sum\limits_{n=1}^{\infty}}
\frac{\varphi(n)}{n^{s+1}}  & =%
{\displaystyle\prod\limits_{p}}
\left(  1+%
{\displaystyle\sum\limits_{r=1}^{\infty}}
\frac{\varphi(p^{r})}{p^{r(s+1)}}\right) \\
& =%
{\displaystyle\prod\limits_{p}}
\left(  1+%
{\displaystyle\sum\limits_{r=1}^{\infty}}
\frac{p-1}{p^{rs+1}}\right) \\
& =%
{\displaystyle\prod\limits_{p}}
\left(  1+\frac{p-1}p%
{\displaystyle\sum\limits_{r=1}^{\infty}}
\frac1{p^{rs}}\right) \\
& =%
{\displaystyle\prod\limits_{p}}
\left(  1+\frac{p-1}{p(p^{s}-1)}\right)  =G(s)\cdot\zeta(s)
\end{align*}
where $G(s)$ is bounded in a half plane $\operatorname*{Re}(s)>c$ for some
$c<1$. In fact, $G(s)=1/\zeta(s+1)$ in this case, and hence
\[%
{\displaystyle\sum\limits_{n\leq N}}
\frac{\varphi(n)}n\sim\frac{G(1)}{\Gamma(1)}N=\frac1{\zeta(2)}N
\]
as $N\rightarrow\infty$. It follows by partial summation that
\[%
{\displaystyle\sum\limits_{n\leq N}}
\varphi(n)\sim\frac1{2\zeta(2)}N^{2}=\frac3{\pi^{2}}N^{2}.
\]
A more elementary proof of this asymptotic formula appears in \cite{Apstl}.
Since the Selberg-Delange method will be used throughout this paper, we choose
to illustrate its application in this simple setting.

Many other questions can be asked for arbitrary $n$:

\begin{itemize}
\item What is the average number of solutions of $x^{2}=1$ in $\mathbb{Z}%
_{n}^{*}?$

\item What is the average number of solutions of $x^{3}=1$ in $\mathbb{Z}%
_{n}^{*}?$

\item What is the average number of solutions of $x^{2}=0$ in $\mathbb{Z}%
_{n}?$

\item What is the average number of solutions of $x^{3}=0$ in $\mathbb{Z}%
_{n}?$

\item What is the average number of images of the map $y\mapsto y^{2}$ in
either $\mathbb{Z}_{n}$ or $\mathbb{Z}_{n}^{*}?$

\item What is the average number of images of the map $y\mapsto y^{3}$ in
either $\mathbb{Z}_{n}$ or $\mathbb{Z}_{n}^{*}?$
\end{itemize}

Although the answers require only straightforward use of standard techniques,
they do not seem to be explicitly given in the literature. We make no claim of
originality:\ Our purpose is only to collect results in one place and to
document relevant numerical techniques.

\section{Number Theory}

\subsection{Selberg-Delange Method}

Let $F(s)=\sum_{n=1}^{\infty}a(n)n^{-s}$ be a Dirichlet series with positive
coefficients and with the property that $G(s)=F(s)\cdot\zeta(s)^{-z}$ can be
analytically continued and is bounded over $\operatorname*{Re}(s)>c$, for some
$c<1$ and some $z\in\mathbb{C}$. Then
\[%
{\displaystyle\sum\limits_{n\leq N}}
a(n)\sim\frac{G(1)}{\Gamma(z)}N\cdot(\ln N)^{z-1}%
\]
as $N\rightarrow\infty$. More terms of the asymptotic expansion are possible,
as is an accurate estimate of the error, but we omit these details for
brevity's sake.

A generalization of this method is required for our work involving averages
over arithmetic progressions. Let $\chi$ denote the principal character modulo
$k=q^{m}$, where $q$ is a prime and $m\geq1$. Here we examine
\[
F_{\chi}(s)=%
{\displaystyle\sum\limits_{n=1}^{\infty}}
\frac{\chi(n)a(n)}{n^{s}}=G_{\chi}(s)\cdot L_{\chi}(s)^{z}%
\]
where
\[
L_{\chi}(s)=%
{\displaystyle\sum\limits_{n=1}^{\infty}}
\frac{\chi(n)}{n^{s}}=\left(  1-\frac1{k^{s}}\right)  \zeta(s)
\]
is the L-series corresponding to $\chi$. Assuming $a(n)$ is a multiplicative
function and $q\nmid\ell$, it follows that
\[%
{\displaystyle\sum\limits_{\substack{n\leq N, \\n\equiv\ell\operatorname*{mod}%
k }}}
a(n)\sim\frac{G_{\chi}(1)}{\Gamma(z)}N\cdot(\ln N)^{z-1}%
\]
as $N\rightarrow\infty$ and, further, that
\[
G_{\chi}(1)=\frac1{\varphi(k)}\left(  1+%
{\displaystyle\sum\limits_{r=1}^{\infty}}
\frac{a(q^{r})}{q^{r}}\right)  ^{-1}G(1).
\]
In our examples, $q$ will be either $2$ or $3$ and the bracketed infinite
series will always collapse to a closed-form expression. Rather than directly
employing the formula for $G_{\chi}(1)$, however, we prefer instead to deduce
$F_{\chi}(s)$ (and hence $G_{\chi}(s)$) from $F(s)$ on basic principles.

\subsection{\label{Sqrt}Square Roots of Unity}

The number $a(n)$ of solutions of $x^{2}=1$ in $\mathbb{Z}_{n}^{*}$ is
\cite{Sec2}
\[
a(n)=\left\{
\begin{array}
[c]{lll}%
2^{\omega(n)-1} &  & \text{if }n\equiv2,6\operatorname*{mod}8,\\
2^{\omega(n)} &  & \text{if }n\equiv1,3,4,5,7\operatorname*{mod}8,\\
2^{\omega(n)+1} &  & \text{if }n\equiv0\operatorname*{mod}8
\end{array}
\right.
\]
where $\omega(n)$ denotes the number of distinct prime factors of $n$. It is
well-known that \cite{Finch}
\begin{align*}%
{\displaystyle\sum\limits_{n=1}^{\infty}}
\frac{2^{\omega(n)}}{n^{s}}  & =%
{\displaystyle\prod\limits_{p}}
\left(  1+2%
{\displaystyle\sum\limits_{r=1}^{\infty}}
\frac1{p^{rs}}\right) \\
\  & =%
{\displaystyle\prod\limits_{p}}
\left(  1+\frac2{p^{s}-1}\right)  =\frac{\zeta(s)^{2}}{\zeta(2s)}%
=G(s)\cdot\zeta(s)^{2}%
\end{align*}
and hence
\[%
{\displaystyle\sum\limits_{n\leq N}}
2^{\omega(n)}\sim\frac1{\zeta(2)}N\cdot\ln N=\frac6{\pi^{2}}N\cdot\ln N.
\]
We need to generalize this asymptotic formula to arithmetic progressions
$n\equiv\ell\operatorname*{mod}k$, where $k=2^{m}$ for $m\geq1$ and
$2\nmid\ell$. It can be shown that
\begin{align*}%
{\displaystyle\sum\limits_{n\equiv\ell\operatorname*{mod}k}}
\frac{2^{\omega(n)}}{n^{s}}  & \sim\frac1{\varphi(k)}%
{\displaystyle\prod\limits_{p>2}}
\left(  1+\frac2{p^{s}-1}\right) \\
\  & =\frac2k\left(  1+\frac2{2^{s}-1}\right)  ^{-1}\frac{\zeta(s)^{2}}%
{\zeta(2s)}=G_{\chi}(s)\cdot\zeta(s)^{2}%
\end{align*}
as $s\rightarrow1$, and thus
\[%
{\displaystyle\sum\limits_{\substack{n\leq N, \\n\equiv\ell\operatorname*{mod}%
k }}}
2^{\omega(n)}\sim\frac{G_{\chi}(1)}{\Gamma(2)}N\cdot\ln N=\frac4{k\pi^{2}%
}N\cdot\ln N.
\]
The cases $(k,\ell)=(8,1),(8,3),(8,5)$ and $(8,7)$ follow immediately. The
case $(k,\ell)=(8,4)$ proceeds from the case $(k,\ell)=(2,1)$:
\[
\frac1{N\cdot\ln N}%
{\displaystyle\sum\limits_{\substack{n\leq N, \\n\equiv4\operatorname*{mod}8
}}}
2^{\omega(n)}=\frac1{N\cdot\ln N}%
{\displaystyle\sum\limits_{\substack{4n\leq N, \\n\equiv1\operatorname*{mod}2
}}}
2^{\omega(n)+1}\longrightarrow\frac14\cdot2\cdot\frac2{\pi^{2}}=\frac1{\pi
^{2}}.
\]
The case $(k,\ell)=(8,2)$ proceeds from the case $(k,\ell)=(4,1)$:
\[
\frac1{N\cdot\ln N}%
{\displaystyle\sum\limits_{\substack{n\leq N, \\n\equiv2\operatorname*{mod}8
}}}
2^{\omega(n)}=\frac1{N\cdot\ln N}%
{\displaystyle\sum\limits_{\substack{2n\leq N, \\n\equiv1\operatorname*{mod}4
}}}
2^{\omega(n)+1}\longrightarrow\frac12\cdot2\cdot\frac1{\pi^{2}}=\frac1{\pi
^{2}}%
\]
and $(8,6)$ likewise proceeds from $(4,1)$. By everything proved thus far, we
have
\[
\frac1{N\cdot\ln N}%
{\displaystyle\sum\limits_{\substack{n\leq N, \\n\equiv0\operatorname*{mod}8
}}}
2^{\omega(n)}\longrightarrow\frac6{\pi^{2}}-4\cdot\frac1{2\pi^{2}}-3\cdot
\frac1{\pi^{2}}=\frac1{\pi^{2}}.
\]
Therefore
\[%
{\displaystyle\sum\limits_{n\leq N}}
a(n)\sim\left(  \frac12\cdot2\cdot\frac1{\pi^{2}}+\left(  4\cdot\frac
1{2\pi^{2}}+\frac1{\pi^{2}}\right)  +2\cdot\frac1{\pi^{2}}\right)  N\cdot\ln
N=\frac6{\pi^{2}}N\cdot\ln N.
\]
It is interesting that $(8,2)$ and $(8,6)$ balance perfectly against $(8,0)$
so that the mean value of $a(n)$ is asymptotically equivalent to the mean
value of $2^{\omega(n)}$.

\subsection{\label{Cbrt}Cube Roots of Unity}

The number $a(n)$ of solutions of $x^{3}=1$ in $\mathbb{Z}_{n}^{*}$ is
\cite{Sec3}
\[
a(n)=\left\{
\begin{array}
[c]{lll}%
3^{\tilde\omega(n)} &  & \text{if }n\equiv1,2,3,4,5,6,7,8\operatorname*{mod}%
9,\\
3^{\tilde\omega(n)+1} &  & \text{if }n\equiv0\operatorname*{mod}9
\end{array}
\right.
\]
where $\tilde\omega(n)$ denotes the number of distinct primes of the form
$3k+1$ dividing $n$:
\[
\tilde\omega(p^{r})=\left\{
\begin{array}
[c]{lll}%
0 &  & \text{if }p=3\text{ or }p\equiv2\operatorname*{mod}3,\\
1 &  & \text{if }p\equiv1\operatorname*{mod}3.
\end{array}
\right.
\]
First, note that
\begin{align*}%
{\displaystyle\sum\limits_{n=1}^{\infty}}
\frac{3^{\tilde\omega(n)}}{n^{s}}  & =%
{\displaystyle\prod\limits_{\substack{p=3\text{ or} \\p\equiv
2\operatorname*{mod}3 }}}
\left(  1+%
{\displaystyle\sum\limits_{r=1}^{\infty}}
\frac1{p^{rs}}\right)  \cdot%
{\displaystyle\prod\limits_{p\equiv1\operatorname*{mod}3}}
\left(  1+3%
{\displaystyle\sum\limits_{r=1}^{\infty}}
\frac1{p^{rs}}\right) \\
\  & =%
{\displaystyle\prod\limits_{\substack{p=3\text{ or} \\p\equiv
2\operatorname*{mod}3 }}}
\left(  1+\frac1{p^{s}-1}\right)  \cdot%
{\displaystyle\prod\limits_{p\equiv1\operatorname*{mod}3}}
\left(  1+\frac3{p^{s}-1}\right) \\
\  & =%
{\displaystyle\prod\limits_{\substack{p=3\text{ or} \\p\equiv
2\operatorname*{mod}3 }}}
\left(  1-\frac1{p^{s}}\right)  ^{-1}\cdot%
{\displaystyle\prod\limits_{p\equiv1\operatorname*{mod}3}}
\left(  1-\frac1{p^{s}}\right)  ^{-3}\left(  1-\frac1{p^{2s}}\right)  \left(
1-\frac2{p^{s}(p^{s}+1)}\right) \\
\  & =\zeta(s)\cdot%
{\displaystyle\prod\limits_{p\equiv1\operatorname*{mod}3}}
\left(  1-\frac1{p^{s}}\right)  ^{-2}\left(  1-\frac1{p^{2s}}\right)  \left(
1-\frac2{p^{s}(p^{s}+1)}\right)  =G(s)\cdot\zeta(s)^{2}%
\end{align*}
and \cite{Moree}
\[
\lim_{s\rightarrow1}%
{\displaystyle\prod\limits_{p\equiv1\operatorname*{mod}3}}
\left(  1-\frac1{p^{s}}\right)  ^{-2}\cdot(s-1)=\frac{\sqrt{3}}{2\pi}%
{\displaystyle\prod\limits_{p\equiv1\operatorname*{mod}3}}
\left(  1-\frac1{p^{2}}\right)  ^{-1};
\]
hence
\[%
{\displaystyle\sum\limits_{n\leq N}}
3^{\tilde\omega(n)}\sim\frac{G(1)}{\Gamma(2)}N\cdot\ln N=C\cdot N\cdot\ln N
\]
where
\[
C=\frac{\sqrt{3}}{2\pi}%
{\displaystyle\prod\limits_{p\equiv1\operatorname*{mod}3}}
\left(  1-\frac2{p(p+1)}\right)  =\frac{\sqrt{3}}{2\pi}%
(0.9410349413195354517900322...).
\]
We need to generalize this asymptotic formula to arithmetic progressions
$n\equiv\ell\operatorname*{mod}k$, where $k=3^{m}$ for $m\geq1$ and
$3\nmid\ell$. It can be shown that
\begin{align*}%
{\displaystyle\sum\limits_{n\equiv\ell\operatorname*{mod}k}}
\frac{3^{\tilde\omega(n)}}{n^{s}}  & \sim\frac1{\varphi(k)}%
{\displaystyle\prod\limits_{p\equiv2\operatorname*{mod}3}}
\left(  1+\frac1{p^{s}-1}\right)  \cdot%
{\displaystyle\prod\limits_{p\equiv1\operatorname*{mod}3}}
\left(  1+\frac3{p^{s}-1}\right) \\
\  & =\frac3{2k}\left(  1+\frac1{3^{s}-1}\right)  ^{-1}G(s)\cdot\zeta
(s)^{2}=G_{\chi}(s)\cdot\zeta(s)^{2}%
\end{align*}
as $s\rightarrow1$, and thus
\[%
{\displaystyle\sum\limits_{\substack{n\leq N, \\n\equiv\ell\operatorname*{mod}%
k }}}
3^{\tilde\omega(n)}\sim\frac{G_{\chi}(1)}{\Gamma(2)}N\cdot\ln N=\frac
CkN\cdot\ln N.
\]
The cases $(k,\ell)=(9,1),(9,2),(9,4),(9,5),(9,7)$ and $(9,8)$ follow
immediately. The case $(9,3)$ proceeds from the case $(3,1)$:
\[
\frac1{N\cdot\ln N}%
{\displaystyle\sum\limits_{\substack{n\leq N, \\n\equiv3\operatorname*{mod}9
}}}
3^{\tilde\omega(n)}=\frac1{N\cdot\ln N}%
{\displaystyle\sum\limits_{\substack{3n\leq N, \\n\equiv1\operatorname*{mod}3
}}}
3^{\tilde\omega(n)}\longrightarrow\frac13\cdot\frac C3=\frac C9
\]
and $(9,6),(9,0)$ likewise proceed from $(3,2),(3,0)$. Therefore
\[%
{\displaystyle\sum\limits_{n\leq N}}
a(n)\sim\left(  8\cdot\frac C9+3\cdot\frac C9\right)  N\cdot\ln N=\frac
{11}9C\cdot N\cdot\ln N=(0.317...)N\cdot\ln N.
\]
Unlike earlier, the mean value of $a(n)$ is asymptotically greater than the
mean value of $3^{\tilde\omega(n)}$. Our estimate improves upon Cloitre
\cite{Sec3}, who gave $(0.4...)N\cdot\ln(N)$ on empirical grounds.

\subsection{Squares in $\mathbb{Z}_{n}^{*}$}

Let $a(n)$ be as defined in section [\ref{Sqrt}]. The number of squares, that
is, the cardinality of images under the map $y\mapsto y^{2}$ in $\mathbb{Z}%
_{n}^{*}$, is \cite{Sec4}
\[
b(n)=\frac{\varphi(n)}{a(n)}=\left\{
\begin{array}
[c]{lll}%
\dfrac{\varphi(n)}{2^{\omega(n)-1}} &  & \text{if }n\equiv
2,6\operatorname*{mod}8,\\
\dfrac{\varphi(n)}{2^{\omega(n)}} &  & \text{if }n\equiv
1,3,4,5,7\operatorname*{mod}8,\\
\dfrac{\varphi(n)}{2^{\omega(n)+1}} &  & \text{if }n\equiv0\operatorname*{mod}%
8.
\end{array}
\right.
\]
First, note that
\begin{align*}%
{\displaystyle\sum\limits_{n=1}^{\infty}}
\frac{\varphi(n)}{n^{s+1}2^{\omega(n)}}  & =%
{\displaystyle\prod\limits_{p}}
\left(  1+\frac12%
{\displaystyle\sum\limits_{r=1}^{\infty}}
\frac{p-1}{p^{rs+1}}\right) \\
\  & =%
{\displaystyle\prod\limits_{p}}
\left(  1+\frac{p-1}{2p(p^{s}-1)}\right)  =G(s)\cdot\zeta(s)^{1/2},
\end{align*}
hence
\[%
{\displaystyle\sum\limits_{n\leq N}}
\frac{\varphi(n)}{n2^{\omega(n)}}\sim\frac{G(1)}{\Gamma(1/2)}N\cdot(\ln
N)^{-1/2}=C\cdot N\cdot(\ln N)^{-1/2}%
\]
where
\[
C=\frac1{\sqrt{\pi}}%
{\displaystyle\prod\limits_{p}}
\left(  1+\frac1{2p}\right)  \left(  1-\frac1p\right)  ^{1/2}=\frac1{\sqrt
{\pi}}(0.8121057111631225117062509...).
\]
It follows by partial summation that
\[%
{\displaystyle\sum\limits_{n\leq N}}
\frac{\varphi(n)}{2^{\omega(n)}}\sim\frac C2\cdot N^{2}\cdot(\ln N)^{-1/2}.
\]
We need to generalize this asymptotic formula to arithmetic progressions
$n\equiv\ell\operatorname*{mod}k$, where $k=2^{m}$ for $m\geq1$ and
$2\nmid\ell$. It can be shown that
\begin{align*}%
{\displaystyle\sum\limits_{n\equiv\ell\operatorname*{mod}k}}
\frac{\varphi(n)}{n^{s+1}2^{\omega(n)}}  & \sim\frac1{\varphi(k)}%
{\displaystyle\prod\limits_{p>2}}
\left(  1+\frac{p-1}{2p(p^{s}-1)}\right) \\
\  & =\frac2k\left(  1+\frac1{4(2^{s}-1)}\right)  ^{-1}G(s)\cdot\zeta
(s)^{1/2}=G_{\chi}(s)\cdot\zeta(s)^{1/2}%
\end{align*}
as $s\rightarrow1$, and thus
\[%
{\displaystyle\sum\limits_{\substack{n\leq N, \\n\equiv\ell\operatorname*{mod}%
k }}}
\frac{\varphi(n)}{n2^{\omega(n)}}\sim\frac{G_{\chi}(1)}{\Gamma(1/2)}N\cdot(\ln
N)^{-1/2}=\frac8{5k}C\cdot N\cdot(\ln N)^{-1/2}%
\]
or
\[%
{\displaystyle\sum\limits_{\substack{n\leq N, \\n\equiv\ell\operatorname*{mod}%
k }}}
\frac{\varphi(n)}{2^{\omega(n)}}\sim\frac4{5k}C\cdot N^{2}\cdot(\ln N)^{-1/2}.
\]
The cases $(k,\ell)=(8,1),(8,3),(8,5)$ and $(8,7)$ follow immediately. The
case $(k,\ell)=(8,4)$ proceeds from the case $(k,\ell)=(2,1)$:
\[
\frac{(\ln N)^{1/2}}{N^{2}}%
{\displaystyle\sum\limits_{\substack{n\leq N, \\n\equiv4\operatorname*{mod}8
}}}
\frac{\varphi(n)}{2^{\omega(n)}}=\frac{(\ln N)^{1/2}}{N^{2}}%
{\displaystyle\sum\limits_{\substack{4n\leq N, \\n\equiv1\operatorname*{mod}2
}}}
\frac{2\varphi(n)}{2^{\omega(n)+1}}\longrightarrow\frac1{16}\cdot\frac
4{10}C=\frac C{40}.
\]
The case $(k,\ell)=(8,2)$ proceeds from the case $(k,\ell)=(4,1)$:
\[
\frac{(\ln N)^{1/2}}{N^{2}}%
{\displaystyle\sum\limits_{\substack{n\leq N, \\n\equiv2\operatorname*{mod}8
}}}
\frac{\varphi(n)}{2^{\omega(n)}}=\frac{(\ln N)^{1/2}}{N^{2}}%
{\displaystyle\sum\limits_{\substack{2n\leq N, \\n\equiv1\operatorname*{mod}4
}}}
\frac{\varphi(n)}{2^{\omega(n)+1}}\longrightarrow\frac14\cdot\frac12\cdot
\frac4{20}C=\frac C{40}%
\]
and $(8,6)$ likewise proceeds from $(4,1)$. By everything proved thus far, we
have
\[
\frac{(\ln N)^{1/2}}{N^{2}}%
{\displaystyle\sum\limits_{\substack{n\leq N, \\n\equiv0\operatorname*{mod}8
}}}
\frac{\varphi(n)}{2^{\omega(n)}}\longrightarrow\frac C2-4\cdot\frac
C{10}-\frac C{40}-2\cdot\frac C{40}=\frac C{40}.
\]
Therefore
\begin{align*}%
{\displaystyle\sum\limits_{n\leq N}}
b(n)  & \sim\left(  2\cdot2\cdot\frac C{40}+\left(  4\cdot\frac C{10}+\frac
C{40}\right)  +\frac12\cdot\frac C{40}\right)  N^{2}\cdot(\ln N)^{-1/2}\\
\  & =\frac{43}{80}C\cdot N^{2}\cdot(\ln N)^{-1/2}=(0.246...)N^{2}\cdot(\ln
N)^{-1/2}.
\end{align*}

\subsection{Cubes in $\mathbb{Z}_{n}^{*}$}

Let $a(n)$ be as defined in section [\ref{Cbrt}]. The number of cubes, that
is, the cardinality of images under the map $y\mapsto y^{3}$ in $\mathbb{Z}%
_{n}^{*}$, is \cite{Sec5}
\[
b(n)=\frac{\varphi(n)}{a(n)}=\left\{
\begin{array}
[c]{lll}%
\dfrac{\varphi(n)}{3^{\tilde\omega(n)}} &  & \text{if }n\equiv
1,2,3,4,5,6,7,8\operatorname*{mod}9,\\
\dfrac{\varphi(n)}{3^{\tilde\omega(n)+1}} &  & \text{if }n\equiv
0\operatorname*{mod}9.
\end{array}
\right.
\]
First, note that
\begin{align*}%
{\displaystyle\sum\limits_{n=1}^{\infty}}
\frac{\varphi(n)}{n^{s+1}3^{\tilde\omega(n)}}  & =%
{\displaystyle\prod\limits_{\substack{p=3\text{ or} \\p\equiv
2\operatorname*{mod}3 }}}
\left(  1+%
{\displaystyle\sum\limits_{r=1}^{\infty}}
\frac{p-1}{p^{rs+1}}\right)  \cdot%
{\displaystyle\prod\limits_{p\equiv1\operatorname*{mod}3}}
\left(  1+\frac13%
{\displaystyle\sum\limits_{r=1}^{\infty}}
\frac{p-1}{p^{rs+1}}\right) \\
\  & =%
{\displaystyle\prod\limits_{\substack{p=3\text{ or} \\p\equiv
2\operatorname*{mod}3 }}}
\left(  1+\frac{p-1}{p(p^{s}-1)}\right)  \cdot%
{\displaystyle\prod\limits_{p\equiv1\operatorname*{mod}3}}
\left(  1+\frac{p-1}{3p(p^{s}-1)}\right)  =G(s)\cdot\zeta(s)^{2/3},
\end{align*}
hence
\[%
{\displaystyle\sum\limits_{n\leq N}}
\frac{\varphi(n)}{n3^{\tilde\omega(n)}}\sim\frac{G(1)}{\Gamma(2/3)}N\cdot(\ln
N)^{-1/3}=C\cdot N\cdot(\ln N)^{-1/3}%
\]
where
\begin{align*}
C  & =\frac1{\Gamma(2/3)}%
{\displaystyle\prod\limits_{\substack{p=3\text{ or} \\p\equiv
2\operatorname*{mod}3 }}}
\left(  1+\frac1p\right)  \left(  1-\frac1p\right)  ^{2/3}\cdot%
{\displaystyle\prod\limits_{p\equiv1\operatorname*{mod}3}}
\left(  1+\frac1{3p}\right)  \left(  1-\frac1p\right)  ^{2/3}\\
\  & =\frac1{\Gamma(2/3)}(0.9477556177621765519078142...).
\end{align*}
It follows by partial summation that
\[%
{\displaystyle\sum\limits_{n\leq N}}
\frac{\varphi(n)}{3^{\tilde\omega(n)}}\sim\frac C2\cdot N^{2}\cdot(\ln
N)^{-1/3}.
\]
We need to generalize this asymptotic formula to arithmetic progressions
$n\equiv\ell\operatorname*{mod}k$, where $k=3^{m}$ for $m\geq1$ and
$3\nmid\ell$. It can be shown that
\begin{align*}%
{\displaystyle\sum\limits_{n\equiv\ell\operatorname*{mod}k}}
\frac{\varphi(n)}{n^{s+1}3^{\tilde\omega(n)}}  & \sim\frac1{\varphi(k)}%
{\displaystyle\prod\limits_{p\equiv2\operatorname*{mod}3}}
\left(  1+\frac{p-1}{p(p^{s}-1)}\right)  \cdot%
{\displaystyle\prod\limits_{p\equiv1\operatorname*{mod}3}}
\left(  1+\frac{p-1}{3p(p^{s}-1)}\right) \\
\  & =\frac3{2k}\left(  1+\frac2{3(3^{s}-1)}\right)  ^{-1}G(s)\cdot
\zeta(s)^{2/3}=G_{\chi}(s)\cdot\zeta(s)^{2/3}%
\end{align*}
as $s\rightarrow1$, and thus
\[%
{\displaystyle\sum\limits_{\substack{n\leq N, \\n\equiv\ell\operatorname*{mod}%
k }}}
\frac{\varphi(n)}{n3^{\tilde\omega(n)}}\sim\frac{G_{\chi}(1)}{\Gamma
(2/3)}N\cdot(\ln N)^{-1/3}=\frac9{8k}C\cdot N\cdot(\ln N)^{-1/3}%
\]
or
\[%
{\displaystyle\sum\limits_{\substack{n\leq N, \\n\equiv\ell\operatorname*{mod}%
k }}}
\frac{\varphi(n)}{3^{\tilde\omega(n)}}\sim\frac9{16k}C\cdot N^{2}\cdot(\ln
N)^{-1/3}.
\]
The cases $(k,\ell)=(9,1),(9,2),(9,4),(9,5),(9,7)$ and $(9,8)$ follow
immediately. The case $(9,3)$ proceeds from the case $(3,1)$:
\[
\frac{(\ln N)^{1/3}}{N^{2}}%
{\displaystyle\sum\limits_{\substack{n\leq N, \\n\equiv3\operatorname*{mod}9
}}}
\frac{\varphi(n)}{3^{\tilde\omega(n)}}=\frac{(\ln N)^{1/3}}{N^{2}}%
{\displaystyle\sum\limits_{\substack{3n\leq N, \\n\equiv1\operatorname*{mod}3
}}}
\frac{2\varphi(n)}{3^{\tilde\omega(n)}}\longrightarrow\frac19\cdot2\cdot
\frac9{48}C=\frac C{24}%
\]
and $(9,6)$ likewise proceeds from $(3,2)$. By everything proved thus far, we
have
\[
\frac{(\ln N)^{1/3}}{N^{2}}%
{\displaystyle\sum\limits_{\substack{n\leq N, \\n\equiv0\operatorname*{mod}9
}}}
\frac{\varphi(n)}{3^{\tilde\omega(n)}}\longrightarrow\frac C2-6\cdot\frac
C{16}-2\cdot\frac C{24}=\frac C{24}.
\]
Therefore
\begin{align*}%
{\displaystyle\sum\limits_{n\leq N}}
b(n)  & \sim\left(  \left(  6\cdot\frac C{16}+2\cdot\frac C{24}\right)
+\frac13\cdot\frac C{24}\right)  N^{2}\cdot(\ln N)^{-1/3}\\
\  & =\frac{17}{36}C\cdot N^{2}\cdot(\ln N)^{-1/3}=(0.330...)N^{2}\cdot(\ln
N)^{-1/3}.
\end{align*}

\subsection{Square Roots of Nullity}

The number $a(n)$ of solutions of $x^{2}=0$ in $\mathbb{Z}_{n}$ is a
multiplicative function of $n$, with $a(p^{r})=p^{\left\lfloor
r/2\right\rfloor }$, thus \cite{Sec6, Brghn}
\begin{align*}%
{\displaystyle\sum\limits_{n=1}^{\infty}}
\frac{a(n)}{n^{s}}  & =%
{\displaystyle\prod\limits_{p}}
\left(  1+\frac1{p^{s}}+\frac p{p^{2s}}+\frac p{p^{3s}}+\frac{p^{2}}{p^{4s}%
}+\frac{p^{2}}{p^{5s}}+\frac{p^{3}}{p^{6s}}+\frac{p^{3}}{p^{7s}}+\cdots\right)
\\
\  & =%
{\displaystyle\prod\limits_{p}}
\left(  1+\frac1{p^{s}}1+\frac p{p^{2s}}+\frac1{p^{s}}\frac p{p^{2s}}%
+\frac{p^{2}}{p^{4s}}+\frac1{p^{s}}\frac{p^{2}}{p^{4s}}+\frac{p^{3}}{p^{6s}%
}+\frac1{p^{s}}\frac{p^{3}}{p^{6s}}+\cdots\right) \\
\  & =%
{\displaystyle\prod\limits_{p}}
\left(  \dfrac1{1-\dfrac p{p^{2s}}}+\dfrac{1/p^{s}}{1-\dfrac p{p^{2s}}%
}\right)  =%
{\displaystyle\prod\limits_{p}}
\left(  1-\dfrac1{p^{2s-1}}\right)  ^{-1}\left(  1+\frac1{p^{s}}\right) \\
\  & =\frac{\zeta(2s-1)\zeta(s)}{\zeta(2s)}=G(s)\cdot\zeta(s)^{2}%
\end{align*}
and $\lim_{s\rightarrow1}\zeta(2s-1)\cdot(s-1)=1/2$, hence
\[%
{\displaystyle\sum\limits_{n\leq N}}
a(n)\sim\frac{G(1)}{\Gamma(2)}N\cdot\ln N=\frac3{\pi^{2}}N\cdot\ln N.
\]

\subsection{Cube Roots of Nullity}

The number $a(n)$ of solutions of $x^{3}=0$ in $\mathbb{Z}_{n}$ is a
multiplicative function of $n$, with $a(p^{r})=p^{\left\lfloor
2r/3\right\rfloor }$, thus \cite{Sec7}
\begin{align*}%
{\displaystyle\sum\limits_{n=1}^{\infty}}
\frac{a(n)}{n^{s}}  & =%
{\displaystyle\prod\limits_{p}}
\left(  1+\frac1{p^{s}}+\frac p{p^{2s}}+\frac{p^{2}}{p^{3s}}+\frac{p^{2}%
}{p^{4s}}+\frac{p^{3}}{p^{5s}}+\frac{p^{4}}{p^{6s}}+\frac{p^{4}}{p^{7s}}%
+\frac{p^{5}}{p^{8s}}+\frac{p^{6}}{p^{9s}}+\frac{p^{6}}{p^{10s}}+\cdots\right)
\\
\  & =%
{\displaystyle\prod\limits_{p}}
\left(  1+\tfrac1{p^{s}}1+\tfrac p{p^{2s}}1+\tfrac{p^{2}}{p^{3s}}%
+\tfrac1{p^{s}}\tfrac{p^{2}}{p^{3s}}+\tfrac p{p^{2s}}\tfrac{p^{2}}{p^{3s}%
}+\tfrac{p^{4}}{p^{6s}}+\tfrac1{p^{s}}\tfrac{p^{4}}{p^{6s}}+\tfrac p{p^{2s}%
}\tfrac{p^{4}}{p^{6s}}+\tfrac{p^{6}}{p^{9s}}+\tfrac1{p^{s}}\tfrac{p^{6}%
}{p^{9s}}+\cdots\right) \\
\  & =%
{\displaystyle\prod\limits_{p}}
\left(  \dfrac1{1-\dfrac{p^{2}}{p^{3s}}}+\dfrac{1/p^{s}}{1-\dfrac{p^{2}%
}{p^{3s}}}+\dfrac{p/p^{2s}}{1-\dfrac{p^{2}}{p^{3s}}}\right)  =%
{\displaystyle\prod\limits_{p}}
\left(  1-\dfrac1{p^{3s-2}}\right)  ^{-1}\left(  1+\frac1{p^{s}}%
+\frac1{p^{2s-1}}\right) \\
\  & =\zeta(3s-2)%
{\displaystyle\prod\limits_{p}}
\left(  1+\frac1{p^{s}}+\frac1{p^{2s-1}}\right)  =G(s)\cdot\zeta(s)^{3}.
\end{align*}
We have
\begin{align*}
\lim_{s\rightarrow1}%
{\displaystyle\prod\limits_{p}}
\left(  1+\frac1{p^{s}}+\frac1{p^{2s-1}}\right)  \cdot(s-1)^{2}  &
=\lim_{s\rightarrow1}\frac1{\zeta(s)\cdot2\zeta(2s-1)}%
{\displaystyle\prod\limits_{p}}
\left(  1+\frac1{p^{s}}+\frac1{p^{2s-1}}\right) \\
\  & =\frac12\lim_{s\rightarrow1}%
{\displaystyle\prod\limits_{p}}
\left(  1+\frac1{p^{s}}+\frac1{p^{2s-1}}\right)  \left(  1-\frac1{p^{s}%
}\right)  \left(  1-\frac1{p^{2s-1}}\right) \\
\  & =\frac12%
{\displaystyle\prod\limits_{p}}
\left(  1+\frac2p\right)  \left(  1-\frac1p\right)  ^{2}\\
\  & =\frac12%
{\displaystyle\prod\limits_{p}}
\left(  1+\frac3{p-1}\right)  \left(  1-\frac1p\right)  ^{3}\\
\  & =\frac12%
{\displaystyle\prod\limits_{p}}
\left(  1-\frac1{p^{2}}\right)  \left(  1-\frac2{p(p+1)}\right) \\
\  & =\frac1{2\zeta(2)}%
{\displaystyle\prod\limits_{p}}
\left(  1-\frac2{p(p+1)}\right)
\end{align*}
and $\lim_{s\rightarrow1}\zeta(3s-2)\cdot(s-1)=1/3$, hence
\[%
{\displaystyle\sum\limits_{n\leq N}}
a(n)\sim\frac{G(1)}{\Gamma(3)}N\cdot(\ln N)^{2}=C\cdot N\cdot(\ln N)^{2}%
\]
where \cite{Finch, Sec7}
\[
C=\frac1{2\pi^{2}}%
{\displaystyle\prod\limits_{p}}
\left(  1-\frac2{p(p+1)}\right)  =\frac1{12}(0.2867474284344787341078927...).
\]

\subsection{Squares in $\mathbb{Z}_{n}$}

The number $b(n)$ of images under the map $y\mapsto y^{2}$ in $\mathbb{Z}_{n}$
is a multiplicative function of $n$, with \cite{Sec8, Prmrs, Stngl}
\[
b(p^{r})=\left\{
\begin{array}
[c]{lll}%
\dfrac13\left(  2^{r-1}+4\right)  &  & \text{if }p=2\text{ and }%
r\equiv0\operatorname*{mod}2,\\
\dfrac13\left(  2^{r-1}+5\right)  &  & \text{if }p=2\text{ and }%
r\equiv1\operatorname*{mod}2,\\
\dfrac1{2(p+1)}\left(  p^{r+1}+p+2\right)  &  & \text{if }p>2\text{ and
}r\equiv0\operatorname*{mod}2,\\
\dfrac1{2(p+1)}\left(  p^{r+1}+2p+1\right)  &  & \text{if }p>2\text{ and
}r\equiv1\operatorname*{mod}2
\end{array}
\right.
\]
and
\begin{align*}
F(s)=%
{\displaystyle\sum\limits_{n=1}^{\infty}}
\frac{b(n)}{n^{s+1}}=\left(  1+%
{\displaystyle\sum\limits_{r=1}^{\infty}}
\frac{b(2^{r})}{2^{r(s+1)}}\right)  \cdot%
{\displaystyle\prod\limits_{p>2}}
\left(  1+%
{\displaystyle\sum\limits_{r=1}^{\infty}}
\frac{b(p^{r})}{p^{r(s+1)}}\right)  .
\end{align*}
The left-hand factor in $F(s)$ simplifies to
\begin{align*}
& \ \ \ 1+\frac13%
{\displaystyle\sum\limits_{i=1}^{\infty}}
\frac{2^{2i-1}+4}{2^{(2i)(s+1)}}+\frac13%
{\displaystyle\sum\limits_{j=1}^{\infty}}
\frac{2^{(2j-1)-1}+5}{2^{(2j-1)(s+1)}}\\
\  & =1+\frac12\left(  \frac{4^{s+1}-3}{(4^{s+1}-1)(4^{s}-1)}+2^{s}%
\frac{2\cdot4^{s+1}-7}{(4^{s+1}-1)(4^{s}-1)}\right) \\
\  & =\left(  1+\frac{2^{2s+1}-2^{s+1}-1}{2^{s+2}(2^{2s+1}-2^{s-1}-1)}\right)
\left(  1-\frac{2^{s+1}+2}{2(2^{s+1}+1)(2^{s+1}-1)}\right)  \left(
1-\frac1{2^{s}}\right)  ^{-1}%
\end{align*}
and the $p^{\text{th}}$ right-hand factor simplifies to
\begin{align*}
& \ \ \ 1+\dfrac1{2(p+1)}%
{\displaystyle\sum\limits_{i=1}^{\infty}}
\dfrac{p^{2i+1}+p+2}{p^{(2i)(s+1)}}+\dfrac1{2(p+1)}%
{\displaystyle\sum\limits_{j=1}^{\infty}}
\dfrac{p^{(2j-1)+1}+2p+1}{p^{(2j-1)(s+1)}}\\
\  & =1+\dfrac1{2(p+1)}\left(  \dfrac{p^{2s+3}+p^{2s+1}+2p^{2s}-2p-2}%
{(p^{2s+2}-1)(p^{2s}-1)}+p^{s+1}\dfrac{p^{2s+2}+2p^{2s+1}+p^{2s}%
-2p-2}{(p^{2s+2}-1)(p^{2s}-1)}\right) \\
\  & =\left(  1-\frac{(p^{s+1}+2)(p-1)}{2(p^{s+1}+1)(p^{s+1}-1)}\right)
\left(  1-\frac1{p^{s}}\right)  ^{-1}.
\end{align*}
We have
\[
F(s)=\zeta(s)\left(  1+\frac{2^{2s+1}-2^{s+1}-1}{2^{s+2}(2^{2s+1}-2^{s-1}%
-1)}\right)
{\displaystyle\prod\limits_{p}}
\left(  1-\frac{(p^{s+1}+2)(p-1)}{2(p^{s+1}+1)(p^{s+1}-1)}\right)
=G(s)\cdot\zeta(s)^{1/2}%
\]
and hence
\[%
{\displaystyle\sum\limits_{n\leq N}}
\frac{b(n)}n\sim\frac{G(1)}{\Gamma(1/2)}N\cdot(\ln N)^{-1/2}=C\cdot N\cdot(\ln
N)^{-1/2}%
\]
where
\begin{align*}
C  & =\frac{17}{16}\frac1{\sqrt{\pi}}%
{\displaystyle\prod\limits_{p}}
\left(  1-\frac{p^{2}+2}{2(p^{2}+1)(p+1)}\right)  \left(  1-\frac1p\right)
^{-1/2}\\
\  & =\frac{17}{16}\frac1{\sqrt{\pi}}(1.2569136102101885959492115...).
\end{align*}
It follows by partial summation that
\[%
{\displaystyle\sum\limits_{n\leq N}}
b(n)\sim\frac C2\cdot N^{2}\cdot(\ln N)^{-1/2}=(0.376...)N^{2}\cdot(\ln
N)^{-1/2}.
\]

\subsection{Cubes in $\mathbb{Z}_{n}$}

The number $b(n)$ of images under the map $y\mapsto y^{3}$ in $\mathbb{Z}_{n}$
is a multiplicative function of $n$, with
\[
b(p^{r})=\left\{
\begin{array}
[c]{lll}%
\dfrac1{13}\left(  3^{r+1}+10\right)  &  & \text{if }p=3\text{ and }%
r\equiv0\operatorname*{mod}3,\\
\dfrac1{13}\left(  3^{r+1}+30\right)  &  & \text{if }p=3\text{ and }%
r\equiv1\operatorname*{mod}3,\\
\dfrac1{13}\left(  3^{r+1}+12\right)  &  & \text{if }p=3\text{ and }%
r\equiv2\operatorname*{mod}3,\\
\dfrac1{p^{2}+p+1}\left(  p^{r+2}+p+1\right)  &  & \text{if }p\equiv
2\operatorname*{mod}3\text{ and }r\equiv0\operatorname*{mod}3,\\
\dfrac1{p^{2}+p+1}\left(  p^{r+2}+p^{2}+p\right)  &  & \text{if }%
p\equiv2\operatorname*{mod}3\text{ and }r\equiv1\operatorname*{mod}3,\\
\dfrac1{p^{2}+p+1}\left(  p^{r+2}+p^{2}+1\right)  &  & \text{if }%
p\equiv2\operatorname*{mod}3\text{ and }r\equiv2\operatorname*{mod}3,\\
\dfrac1{3(p^{2}+p+1)}\left(  p^{r+2}+2p^{2}+3p+3\right)  &  & \text{if
}p\equiv1\operatorname*{mod}3\text{ and }r\equiv0\operatorname*{mod}3,\\
\dfrac1{3(p^{2}+p+1)}\left(  p^{r+2}+3p^{2}+3p+2\right)  &  & \text{if
}p\equiv1\operatorname*{mod}3\text{ and }r\equiv1\operatorname*{mod}3,\\
\dfrac1{3(p^{2}+p+1)}\left(  p^{r+2}+3p^{2}+2p+3\right)  &  & \text{if
}p\equiv1\operatorname*{mod}3\text{ and }r\equiv2\operatorname*{mod}3
\end{array}
\right.
\]
and
\begin{align*}
F(s)=%
{\displaystyle\sum\limits_{n=1}^{\infty}}
\frac{b(n)}{n^{s+1}}=\left(  1+%
{\displaystyle\sum\limits_{r=1}^{\infty}}
\frac{b(3^{r})}{3^{r(s+1)}}\right)  \cdot%
{\displaystyle\prod\limits_{\substack{\text{ }p\equiv2 \\\operatorname*{mod}3
}}}
\left(  1+%
{\displaystyle\sum\limits_{r=1}^{\infty}}
\frac{b(p^{r})}{p^{r(s+1)}}\right)  \cdot%
{\displaystyle\prod\limits_{\substack{p\equiv1 \\\operatorname*{mod}3 }}}
\left(  1+%
{\displaystyle\sum\limits_{r=1}^{\infty}}
\frac{b(p^{r})}{p^{r(s+1)}}\right)  .
\end{align*}
The expressions for $b(p^{r})$ follow from a conjecture by Wilson \cite{Sec9};
a proof for the case $p=2$, $r\equiv0\operatorname*{mod}3$ was given by Wilmer
\&\ Schirokauer \cite{WS}. The left-hand factor in $F(s)$ simplifies to
\begin{align*}
& 1+\frac1{13}\left(
{\displaystyle\sum\limits_{i=1}^{\infty}}
\frac{3^{3i+1}+10}{3^{(3i)(s+1)}}+%
{\displaystyle\sum\limits_{j=1}^{\infty}}
\frac{3^{(3j-2)+1}+30}{3^{(3j-2)(s+1)}}+%
{\displaystyle\sum\limits_{k=1}^{\infty}}
\frac{3^{(3k-1)+1}+12}{3^{(3k-1)(s+1)}}\right) \\
& =1+\left(  \frac{7\cdot27^{s}-1}{(27^{s+1}-1)(27^{s}-1)}+3^{2s+1}%
\frac{9\cdot27^{s}-7}{(27^{s+1}-1)(27^{s}-1)}+3^{s+1}\frac{3\cdot27^{s}%
-1}{(27^{s+1}-1)(27^{s}-1)}\right) \\
& =\left(  1-\frac{2\left(  3^{s+2}+1\right)  }{(3^{s+1}+3^{(s+1)/2}%
+1)(3^{s+1}-3^{(s+1)/2}+1)(3^{s+1}-1)}\right)  \left(  1-\frac1{3^{s}}\right)
^{-1}.
\end{align*}
The $p^{\text{th}}$ right-hand factor simplifies to
\begin{align*}
& 1+\dfrac1{p^{2}+p+1}\left(
{\displaystyle\sum\limits_{i=1}^{\infty}}
\dfrac{p^{3i+2}+p+1}{p^{(3i)(s+1)}}+%
{\displaystyle\sum\limits_{j=1}^{\infty}}
\dfrac{p^{(3j-2)+2}+p^{2}+p}{p^{(3j-2)(s+1)}}+%
{\displaystyle\sum\limits_{k=1}^{\infty}}
\dfrac{p^{(3k-1)+2}+p^{2}+1}{p^{(3k-1)(s+1)}}\right) \\
& =1+\dfrac1{p^{2}+p+1}\left(  \dfrac{p^{3s+5}+p^{3s+1}+p^{3s}-p^{2}%
-p-1}{(p^{3s+3}-1)(p^{3s}-1)}\right. \\
& \left.  \;\;+\,p^{2s+2}\dfrac{p^{3s+3}+p^{3s+2}+p^{3s+1}-p^{2}%
-p-1}{(p^{3s+3}-1)(p^{3s}-1)}\right. \\
& \left.  \;\;+\,p^{s+1}\dfrac{p^{3s+4}+p^{3s+2}+p^{3s}-p^{2}-p-1}%
{(p^{3s+3}-1)(p^{3s}-1)}\right) \\
& =\left(  1-\frac{(p^{s+1}+1)(p-1)}{(p^{s+1}+p^{(s+1)/2}+1)(p^{s+1}%
-p^{(s+1)/2}+1)(p^{s+1}-1)}\right)  \left(  1-\frac1{p^{s}}\right)  ^{-1}%
\end{align*}
when $p\equiv2\operatorname*{mod}3$ and
\begin{align*}
& \ \ 1+\tfrac1{3(p^{2}+p+1)}\left(
{\displaystyle\sum\limits_{i=1}^{\infty}}
\tfrac{p^{3i+2}+2p^{2}+3p+3}{p^{(3i)(s+1)}}+%
{\displaystyle\sum\limits_{j=1}^{\infty}}
\tfrac{p^{(3j-2)+2}+3p^{2}+3p+2}{p^{(3j-2)(s+1)}}+%
{\displaystyle\sum\limits_{k=1}^{\infty}}
\tfrac{p^{(3k-1)+2}+3p^{2}+2p+3}{p^{(3k-1)(s+1)}}\right) \\
& =1+\dfrac1{3(p^{2}+p+1)}\left(  \dfrac{p^{3s+5}+2p^{3s+2}+3p^{3s+1}%
+3p^{3s}-3p^{2}-3p-3}{(p^{3s+3}-1)(p^{3s}-1)}\right. \\
& \left.  \;\;+\,p^{2s+2}\dfrac{p^{3s+3}+3p^{3s+2}+3p^{3s+1}+2p^{3s}%
-3p^{2}-3p-3}{(p^{3s+3}-1)(p^{3s}-1)}\right. \\
& \left.  \;\;+\,p^{s+1}\dfrac{p^{3s+4}+3p^{3s+2}+2p^{3s+1}+3p^{3s}%
-3p^{2}-3p-3}{(p^{3s+3}-1)(p^{3s}-1)}\right) \\
& =\left(  1-\frac{(2p^{2s+2}+3p^{s+1}+3)(p-1)}{3(p^{s+1}+p^{(s+1)/2}%
+1)(p^{s+1}-p^{(s+1)/2}+1)(p^{s+1}-1)}\right)  \left(  1-\frac1{p^{s}}\right)
^{-1}%
\end{align*}
when $p\equiv1\operatorname*{mod}3$. We have
\begin{align*}
F(s)  & =\zeta(s)\left(  1-\frac{2\left(  3^{s+2}+1\right)  }{(3^{s+1}%
+3^{(s+1)/2}+1)(3^{s+1}-3^{(s+1)/2}+1)(3^{s+1}-1)}\right) \\
& \;\cdot%
{\displaystyle\prod\limits_{\substack{\text{ }p\equiv2 \\\operatorname*{mod}3
}}}
\left(  1-\frac{(p^{s+1}+1)(p-1)}{(p^{s+1}+p^{(s+1)/2}+1)(p^{s+1}%
-p^{(s+1)/2}+1)(p^{s+1}-1)}\right) \\
& \;\cdot%
{\displaystyle\prod\limits_{\substack{p\equiv1 \\\operatorname*{mod}3 }}}
\left(  1-\frac{(2p^{2s+2}+3p^{s+1}+3)(p-1)}{3(p^{s+1}+p^{(s+1)/2}%
+1)(p^{s+1}-p^{(s+1)/2}+1)(p^{s+1}-1)}\right) \\
& =G(s)\cdot\zeta(s)^{2/3}%
\end{align*}
and hence
\[%
{\displaystyle\sum\limits_{n\leq N}}
\frac{b(n)}n\sim\frac{G(1)}{\Gamma(2/3)}N\cdot(\ln N)^{-1/3}=C\cdot N\cdot(\ln
N)^{-1/3}%
\]
where
\begin{align*}
C  & =\frac{12}{13}\frac1{\Gamma(2/3)}\left(  1-\frac13\right)  ^{-1/3}\\
& \;\;\cdot%
{\displaystyle\prod\limits_{\substack{\text{ }p\equiv2 \\\operatorname*{mod}3
}}}
\left(  1-\frac{p^{2}+1}{(p^{2}+p+1)(p^{2}-p+1)(p+1)}\right)  \left(
1-\frac1p\right)  ^{-1/3}\\
& \;\;\cdot%
{\displaystyle\prod\limits_{\substack{p\equiv1 \\\operatorname*{mod}3 }}}
\left(  1-\frac{2p^{4}+3p^{2}+3}{3(p^{2}+p+1)(p^{2}-p+1)(p+1)}\right)  \left(
1-\frac1p\right)  ^{-1/3}\\
& =\frac{12}{13}\frac1{\Gamma(2/3)}(1.4225831466986636811460982...).
\end{align*}
It follows by partial summation that
\[%
{\displaystyle\sum\limits_{n\leq N}}
b(n)\sim\frac C2\cdot N^{2}\cdot(\ln N)^{-1/3}=(0.484...)N^{2}\cdot(\ln
N)^{-1/3}.
\]
We emphasize that this result is only conjectural.

\subsection{Other Problems}

The power of the Selberg-Delange method is evident (many deeper applications
occur elsewhere in the literature). We merely mention that the number $a(n)$
of solutions of $x^{2}=-1$ in $\mathbb{Z}_{n}^{*}$ satisfies
\[%
{\displaystyle\sum\limits_{n\leq N}}
a(n)\sim\frac3{2\pi}N;
\]
in particular, $x^{2}=-1$ has asymptotically far fewer solutions than
$x^{2}=1$. Such asymmetry does not occur for $x^{3}=\pm1$ (just replace $x$ by
$-x$). See other modular polynomial equations at \cite{Sec10} and the
enumeration of weakly primitive Dirichlet characters at \cite{Fnch, Sec11}.

A more difficult exercise concerns the number $b(n)$ of elements of
$\mathbb{Z}_{n}$ that are \textit{both} squares and cubes. If $w=z^{6}$, then
clearly $w=(z^{3})^{2}=(z^{2})^{3}$. Conversely, if $w=u^{2}=v^{3}$, then
$(uv^{-1})^{6}=(u^{2})^{3}(v^{3})^{-2}=w^{3}w^{-2}=w$. Hence $b(n)$ is the
same as the number of sixth-powers in $\mathbb{Z}_{n}$. Wilson's conjecture
again provides expressions for $b(p^{r})$, which in turn give formulas for
$F(s)$ and $G(s)$. The details of this and other higher-power problems are
left to someone else \cite{Sec12}.\newpage\ \ 

\section{Numerical Techniques}

\subsection{Prime Products}

Here is a method for evaluating constants of the form
\[
C=\prod_{p\equiv\ell\operatorname*{mod}k}f(p)
\]
to high precision, where the product is taken over all primes of the form
$p=mk+\ell$. Suppose that the function $\ln f$ has asymptotic expansion
\[
\ln f(p)=\frac{c_{2}}{p^{s_{2}}}+\frac{c_{3}}{p^{s_{3}}}+\cdots+\frac{c_{n}%
}{p^{s_{n}}}+\cdots
\]
as $p\rightarrow\infty,$ where $\left(  c_{n},s_{n}\right)  $ are real numbers
and $1<s_{2}<...<s_{n}<...$. (Often $s_{n}=n$ occurs.) Define the
$(k,\ell)^{\text{th}}$ \textit{prime zeta function}
\[
P_{k,\ell}(s)=\sum_{p\equiv\ell\operatorname*{mod}k}\frac1{p^{s}}%
\]
for $\operatorname*{Re}(s)>1$; it follows that
\[
\ln C=\sum_{p\equiv\ell\operatorname*{mod}k}\ln f(p)=\sum_{n\geq2}%
c_{n}P_{k,\ell}(s_{n}).
\]
Let $p_{k,\ell}$ denote the smallest prime of the form $mk+\ell$; clearly
$P_{k,\ell}(n)\sim1/p_{k,\ell}^{n}$ as $n\rightarrow\infty$. Consequently, if
the coefficients $c_{n}$ are uniformly bounded, the convergence of the sum is
fast (geometric). It hence remains to accurately compute the values
$P_{k,\ell}(s_{n})$.

\subsection{Prime Zeta Functions}

Let $\operatorname*{Re}(s)>1$. The classical prime zeta function
$P(s)=P_{1,0}(s)$ can be related to the classical zeta function by Euler's
famous product:
\[
\ln\zeta(s)=-\sum_{p}\ln\left(  1-\frac1{p^{s}}\right)  =\sum_{p}\sum_{n\geq
1}\frac1{np^{ns}}=\sum_{n\geq1}\frac{P(ns)}n.
\]
Applying the M\"obius inversion formula, we obtain \cite{Frbrg, Cohen}
\[
P(s)=\sum_{n=1}^{\infty}\frac{\mu(n)}n\ln\zeta(ns).
\]
Since $\ln\zeta(ns)\sim2^{-ns}$ as $n\rightarrow\infty$, only a few terms in
this series are required to compute an accurate value of $P(s).$ Also
$P(s)\sim-\ln(s-1)$ as $s\rightarrow1^{+}$. These facts are useful in
computing constants of the form $\prod_{p}f(p)$.

For constants of the form $\prod_{p\equiv\ell\operatorname*{mod}3}f(p)$, we
need $P_{3,1}(s)$ and $P_{3,2}(s)$. To achieve this, it is necessary to
introduce the two characters modulo $3$:
\[%
\begin{array}
[c]{ccc}%
\chi_{0}(n)=\left\{
\begin{array}
[c]{ccc}%
1 &  & \text{if }n\equiv1\operatorname*{mod}3,\\
1 &  & \text{if }n\equiv2\operatorname*{mod}3,\\
0 &  & \text{if }n\equiv0\operatorname*{mod}3,
\end{array}
\right.  &  & \chi_{1}(n)=\left\{
\begin{array}
[c]{ccc}%
1 &  & \text{if }n\equiv1\operatorname*{mod}3,\\
-1 &  & \text{if }n\equiv2\operatorname*{mod}3,\\
0 &  & \text{if }n\equiv0\operatorname*{mod}3,
\end{array}
\right.
\end{array}
\]
and their associated Dirichlet L-series:
\[%
\begin{array}
[c]{ccc}%
L_{j}(s)=L_{\chi_{j}}(s)=\sum\limits_{n=1}^{\infty}\dfrac{\chi_{j}(n)}{n^{s}%
}=\dfrac{1}{3^{s}}\left(  \chi_{j}(1)\zeta\left(  s,\frac{1}{3}\right)
+\chi_{j}(2)\zeta\left(  s,\frac{2}{3}\right)  \right)  , &  & j=0,1
\end{array}
\]
where $\zeta\left(  s,a\right)  $ is the Hurwitz zeta-function. By well-known
acceleration procedures, series of this nature can be evaluated to many
decimal places.

From the Euler product expressions
\[
L_{0}(s)=\prod_{p}\left(  1-\frac{\chi_{0}(p)}{p^{s}}\right)  ^{-1}%
=\prod_{p\equiv1\operatorname*{mod}3}\left(  1-\frac1{p^{s}}\right)
^{-1}\prod_{p\equiv2\operatorname*{mod}3}\left(  1-\frac1{p^{s}}\right)
^{-1},
\]
\[
L_{1}(s)=\prod_{p}\left(  1-\frac{\chi_{1}(p)}{p^{s}}\right)  ^{-1}%
=\prod_{p\equiv1\operatorname*{mod}3}\left(  1-\frac1{p^{s}}\right)
^{-1}\prod_{p\equiv2\operatorname*{mod}3}\left(  1+\frac1{p^{s}}\right)
^{-1},
\]
we obtain
\[
\frac12\ln\left(  \frac{L_{0}(s)}{L_{1}(s)}\right)  =\frac12\sum
_{p\equiv2\operatorname*{mod}3}\ln\left(  \frac{1+p^{-s}}{1-p^{-s}}\right)
=\sum_{n=0}^{\infty}\frac{P_{3,2}((2n+1)s)}{2n+1},
\]
\[
\frac12\ln\left(  \frac{L_{0}(s)L_{1}(s)}{L_{0}(2s)}\right)  =\frac
12\sum_{p\equiv1\operatorname*{mod}3}\ln\left(  \frac{1+p^{-s}}{1-p^{-s}%
}\right)  =\sum_{n=0}^{\infty}\frac{P_{3,1}((2n+1)s)}{2n+1}%
\]
and, again, by M\"obius inversion,
\[
P_{3,2}(s)=\frac12\sum_{n=0}^{\infty}\frac{\mu(2n+1)}{2n+1}\ln\left(
\frac{L_{0}((2n+1)s)}{L_{1}((2n+1)s)}\right)  ,
\]
\[
P_{3,1}(s)=\frac12\sum_{n=0}^{\infty}\frac{\mu(2n+1)}{2n+1}\ln\left(
\frac{L_{0}((2n+1)s)L_{1}((2n+1)s)}{L_{0}((4n+2)s)}\right)  .
\]
Of course, $L_{0}(s)=\zeta(s)(1-1/3^{s})$ and $L_{1}(s)=1-1/2^{s}%
+1/4^{s}-1/5^{s}+\cdots.$ Also, $P_{3,1}(s)\sim-\frac12\ln(s-1)$ and
$P_{3,2}(s)\sim-\frac12\ln(s-1)$ as $s\rightarrow1^{+}$.

Similar techniques involving characters modulo $k$ can be used to compute
constants of the form $\prod_{p\equiv\ell\operatorname*{mod}k}f(p)$, but for
brevity's sake we do not discuss these here.

\subsection{A Simple Example}

Let us compute the constant
\[
C=\prod_{p\equiv1\operatorname*{mod}3}\left(  1-\frac2{p(p+1)}\right)
=\prod_{p\equiv1\operatorname*{mod}3}\left(  \frac{(p-1)(p+2)}{(p+1)p}\right)
\]
that appears in section [\ref{Cbrt}]. It is easy to establish that
\[
\ln C=\sum_{p\equiv1\operatorname*{mod}3}\left(  \ln\left(  \frac{p-1}%
{p+1}\right)  +\ln\left(  1+\frac2p\right)  \right)  =\sum_{n\geq2}\frac
{c_{n}}nP_{3,1}(n)
\]
where $c_{n}=2^{n}-2$ when $n$ is odd and $c_{n}=-2^{n}$ when $n$ is even.
Since $c_{n}=O(2^{n})$ and $P_{3,1}(n)=O(7^{-n})$, it is more efficient to
compute directly the product up to a certain cutoff $p_{c}$. For example, if
we take $p_{c}=31$, we find
\[
C=\frac{3247695}{3430336}\prod_{\substack{p\equiv1\operatorname*{mod}3, \\p>31
}}\left(  1-\frac2{p(p+1)}\right)
\]
and consequently
\[
\ln C=\ln\frac{3247695}{3430336}+\sum_{n=2}^{\infty}\frac{c_{n}}n\left(
P_{3,1}(n)-\frac1{7^{n}}-\frac1{13^{n}}-\frac1{19^{n}}-\frac1{31^{n}}\right)
\]
enjoys much faster convergence because $P_{3,1}(n)-\frac1{7^{n}}-\frac
1{13^{n}}-\frac1{19^{n}}-\frac1{31^{n}}=O(37^{-n})$. The first few terms of
this series produce
\[%
\begin{tabular}
[c]{c|c}%
$n$ & \multicolumn{1}{c}{$C$}\\\hline
\multicolumn{1}{l|}{$2$} &
\multicolumn{1}{|l}{$0.94(09438379523896292195206...)$}\\
\multicolumn{1}{l|}{$3$} &
\multicolumn{1}{|l}{$0.94103(87732177050567463275...)$}\\
\multicolumn{1}{l|}{$4$} &
\multicolumn{1}{|l}{$0.941034(8096648041499806620...)$}\\
\multicolumn{1}{l|}{$5$} &
\multicolumn{1}{|l}{$0.94103494(70255355752383278...)$}\\
\multicolumn{1}{l|}{$10$} &
\multicolumn{1}{|l}{$0.94103494131953(43277214763...)$}\\
\multicolumn{1}{l|}{$15$} &
\multicolumn{1}{|l}{$0.941034941319535451790(3566...)$}%
\end{tabular}
\]
and only 15 terms are necessary to obtain 20 correct decimal places.

\section{Acknowledgements}

We thank G\'{e}rald Tenenbaum \& Jean-Marie De Koninck for their expertise in
the Selberg-Delange method, Pieter Moree for his help in evaluating a complex
residue, and David Wilson \&\ Benoit Cloitre for their many contributions to
Neil Sloane's sequence database. \ Our work is extended in \cite{FMS} and we
gratefully acknowledge Greg Martin for his mastery of the subject.

\end{document}